\documentclass{article}
\usepackage{graphicx}
\usepackage{amsmath}
\usepackage{amsthm}
\usepackage{amsfonts}
\usepackage{amssymb}
\newtheorem*{theorem}{Theorem}
\newtheorem*{theorem1}{Theorem [11]}
\newtheorem*{corollary}{Corollary}
\newtheorem*{definition}{Definition}

\newtheorem*{lemma}{Lemma}

\newtheorem*{proposition}{Proposition}
\newtheorem*{remark}{Remark}

\begin{document}

\title{Simple completable contractions of nilpotent Lie algebras}
\author{Rutwig Campoamor-Stursberg\thanks{e-mail: rutwig@nfssrv.mat.ucm.es} \\
Departamento de Geometr\'{\i}a y Topolog\'{\i}a\\
Fac. CC. Matem\'{a}ticas Univ. Complutense\\
28040 Madrid ( Spain )}
\date{}
\maketitle

\begin{abstract}
We study a certain class of non-maximal rank contractions of the nilpotent Lie algebra $\frak{g}_{m}$ and show that these contractions are completable Lie algebras. As a consequence a family of solvable complete Lie algebras of non-maximal rank is given in arbitrary dimension.
.\newline AMS Math. Subj. Class. 17B10, 17B30.
\newline \textit{Keywords : contraction, complete Lie algebra}
\end{abstract}

\section{Generalities}

The notion of contraction of a Lie algebra ( also called degeneration by some authors) was originally introduced by physicists ( Segal) as a tool to relate classical and quatum mechanics. Inonu and Wigner used contractions attending to a particularization, namely, that a subalgebra remains fixed through the contraction. This concept, quite restrictive for some purposes, was later generalized by Saletan [13] and Levi-Nahas. The relation between contractions and deformation theory is an important, but not fully exploited question [7]. Orbit closures of Lie algebras, where contractions play a central role, is not an easy problem, and it constitutes an essential tool in the analysis of the componets of the varieties $\frak{L}^{n}$ and $\frak{N}^{n}$.

\bigskip

Let $\frak{L}^{n}$ be the set of complex  Lie algebra laws in dimension $n$.
We identify each law with its structure constants $C_{ij}^{k}$ on a fixed
basis $\left\{  X_{i}\right\}  $ of $\mathbb{C}^{n}$. The Jacobi identities%
\[
\sum_{l=1}^{n}C_{ij}^{k}C_{kl}^{s}+C_{jk}^{l}C_{il}^{s}+C_{ki}^{l}C_{jl}^{s}=0
\]
for $1\leq i\leq j<k\leq n,\;1\leq s\leq n$ show that $\frak{L}^{n}$ is an
algebraic variety. The nilpotent Lie algebra laws $\frak{N}^{n}$ are a closed
subset in $\frak{L}^{n}$. The linear group $GL\left(  n,\mathbb{C}\right)  $
acts on $\frak{L}^{n}$ via changes of basis, i.e., $\left(g*\mu\right)\left(x,y\right)=g\left(\mu\left(g^{-1}\left(x\right),g^{-1}\left(y\right)\right)\right)$ for $g\in GL\left(n,\mathbb{C}\right)$. Let $\mathcal{O}\left(  \mu\right)$ denote the orbit of the law $\mu$ by this action, consisting of all structures in a single isomorphism class.

\begin{definition}
A Lie algebra $\lambda$ is a contraction of a law $\mu$ if $\lambda
\in\overline{O\left(  \mu\right)  }$. 
\end{definition}

Here the topology on the variety is either the metric topology or the Zariski topology. Both topologies lead to the same contractions. A contraction will be denoted by $\mu\rightarrowtail\lambda$. It follows that
the entire orbit $O\left(  \lambda\right)  $ lies in the Zariski-closure of $O\left(
\mu\right)$. In particular, the following condition implies  $\mu
\rightarrowtail\lambda$:%
\[
\exists g_{t}\in GL\left(  n,\mathbb{C}\left(  t\right)  \right)
\text{\ \ such that \ \ }\lim_{t\rightarrow0}\,g_{t}\ast\mu=\lambda
\]
Contractions are also transitive, i.e., if $\lambda\rightarrowtail \mu$ and $\mu\rightarrowtail \psi$, then $\lambda\rightarrowtail \psi$. Thus not any existing contraction must be shown directly.

For a Lie algebra $\lambda$ to be a contraction of $\mu$, the following
conditions are necessary:

\begin{lemma}
Let $\mu\rightarrowtail\lambda$. Then

\begin{enumerate}
\item $\dim\;Der\left(  \mu\right)  <\dim\;Der\left(  \lambda\right)  $

\item $\dim\;\left[  \lambda,\lambda\right]  \leq\dim\;\left[  \mu,\mu\right]  $

\item $\dim\;Z\left(  \lambda\right)  \geq\dim\;Z\left(  \mu\right)  $

\item $rank\left(  \lambda\right)  \geq rank\left(  \mu\right)  $
\end{enumerate}
where $Der\left(\mu\right)$ denotes the algebra of derivations, $Z\left(\mu\right)$ the center of the algebra $\mu$ and $rank\left(\mu\right)$ is the dimension of a maximal toral subalgebra of $\mu$. 
\end{lemma}

Proofs of the assertions can be found in [5], [7] and [14].\newline For a nilpotent Lie algebra $\frak{g}$, we denote the ideals of the central descending sequence by $C^{i}\left(\frak{g}\right)$, i.e., $C^{i}\left(\frak{g}\right)=[\frak{g},C^{i-1}\left(\frak{g}\right)]$ for $i\geq 1$ and $C^{0}\left(\frak{g}\right)=\frak{g}$. If the algebra $\frak{g}$ contracts to an algebra $\frak{h}$, it follows from the lemma above that $dim\left(C^{i}\left(\frak{g}\right)\right)\geq dim \left(C^{i}\left(\frak{g}\right)\right)$. Therefore, if $\frak{g}$ has nilindex $p$, then the nilindex of $\frak{h}$ is $\leq p$.\newline Although we will not make explicit use of the characteristic sequence, it is convenient to recall this invariant: Consider a complex nilpotent Lie algebra $\frak{g}=\left(  \mathbb{C}^{n}%
,\mu\right)  $. For each $X\in\mathbb{C}^{n}$ we denote $c\left(  X\right)  $
the ordered sequence of dimensions of Jordan blocks of the adjoint operator
$ad_{\mu}\left(  X\right)  $.

\begin{definition}
The characteristic sequence of $\frak{g}$ is an isomorphism invariant
$c\left(  \frak{g}\right)  $ defined as
\[
c\left(  \frak{g}\right)  =\sup_{X\in\frak{g}-C^{1}\frak{g}}\left\{  c\left(
X\right)  \right\}
\]
where $C^{1}\frak{g}$ denotes the derived subalgebra.\newline A characteristic sequence is called linear if there exists an integer $n$ such that $c\left(\frak{g}\right)=\left(n,1,..,1\right)$.
\end{definition}

\section{The algebras $\frak{g}_{m}\left(q_{1},..,q_{k}\right)$}
In this section we analyze some families of nilpotent Lie algebras for which certain classes of contractions will be determined.
\bigskip

For $m\geq4$ let $\frak{g}_{m}$ be the Lie algebra
whose structural equations are
\begin{align*}
d\omega_{1}  &  =d\omega_{2}=0\\
d\omega_{j}  &  =\omega_{1}\wedge\omega_{j-1},\;3\leq j\leq2m\\
d\omega_{2m+1}  &  =\sum_{j=2}^{m}\left(  -1\right)  ^{j}\omega_{j}%
\wedge\omega_{2m+1-j}%
\end{align*}
where $\left\{  \omega_{1},..,\omega_{2m+1}\right\}  $ is a basis of $\left(
\mathbb{C}^{2m+1}\right)  ^{\ast}$.

\begin{proposition}
For any $m\geq4$ the Lie algebra $\frak{g}_{m}$ is
naturally graded of characteristic sequence $\left(  2m-1,1,1\right)$ satisfying the following property
\begin{align*}
C_{\frak{g}_{m}}\left(C^{m}\left(\frak{g}_{m}\right)\right)\supset C^{m}\left(\frak{g}_{m}\right)\\ 
C_{\frak{g}_{m}}\left(C^{m-1}\left(\frak{g}_{m}\right)\right)\not\supset C^{m-1}\left(\frak{g}_{m}\right)
\end{align*}
\end{proposition}

This algebra, which has been analyzed in [3], can be characterized as follows: 

\begin{theorem}
For $m\geq4$ any naturally graded, central extension $\frak{g}$ of the filiform model Lie algebra $L_{2m}$ whose nilindex is $\left(  2m-1\right)$ and satisfies
\begin{align*}
C_{\frak{g}}\left(C^{m}\left(\frak{g}\right)\right)\supset C^{m}\left(\frak{g}\right)\\
C_{\frak{g}}\left(C^{m-1}\left(\frak{g}\right)\right)\not\supset C^{m-1}\left(\frak{g}\right)
\end{align*}
is isomorphic to $\frak{g}_{m}$.
\end{theorem}

Algebras satisfying the preceding "centralizer condition" arise from the analysis of gradations of nilradicals of Borel subalgebras of complex simple Lie algebras, and are defined by a modification of the graded structure of these algebras [3]. In particular, the deformations and extensions of $\frak{g}_{m}$ have been studied in [4].\newline

Let $m\geq4$. For any sequence $3\leq q_{1}<q_{2}<...<q_{k}\leq m+1$ let
$\frak{g}_{m}\left(  q_{1},.,q_{k}\right)  $ be the $\left(  2m+1\right)
$-dimensional Lie algebra whose Maurer-Cartan equations are:%
\begin{align*}
d\omega_{1}  & =d\omega_{2}=0\\
d\omega_{q_{i}}  & =d\omega_{2m+2-q_{i}}=0,\;1\leq i\leq k\\
d\omega_{j}  & =\omega_{1}\wedge\omega_{j-1},\;3\leq j\leq2m,\;j\notin
\{q_{i},2m+2-q_{i}\}_{1\leq i\leq k}\\
d\omega_{2m+1}  & =\sum_{j=2}^{m}\left(  -1\right)  ^{j}\omega_{j}\wedge
\omega_{2m+1-j}%
\end{align*}
where $\left\{  \omega_{1},..,\omega_{2m+1}\right\}  $ is a basis of $\left(
\mathbb{C}^{2m+1}\right)  ^{\ast}$.

\begin{lemma}
For $m\geq4$ and $k\geq1$ the Lie algebras $\frak{g}_{m}\left(  q_{1}%
,.,q_{k}\right)  $ are nonsplit nilpotent of nonlinear characteristic sequence.
\end{lemma}

\begin{remark}
In general, the algebras $\frak{g}_{m}\left(  q_{1},.,q_{k}\right)  $ will not
be naturally graded. In particular this happens whenever we have $q,q^{\prime
}$ such that $q^{\prime}=1+q$. In fact, the differential form $d\omega_{2m+1}$ determines the gradation in some sense [3]. These algebras are also interesting for the study of solvable rigid Lie algebras whose nilradical has linear characteristic sequence [2].
\end{remark}

\begin{theorem}
For any $m\geq4$ and $k\geq1$ the Lie algebra $\frak{g}_{m}\left(
q_{1},.,q_{k}\right)  $ is a contraction of $\frak{g}$.
\end{theorem}

\begin{proof}
For any $k$-tuple $\left(  q_{1},..,q_{k}\right)  $ let $S\left(
m,q_{1},..,q_{k}\right)  $ be the following linear system
\begin{align}
a_{1}+a_{j-1}  & =a_{j},\;3\leq j\leq2m,\;j\notin\left\{  q_{i},2m+2-q_{i}%
\right\}  _{1\leq i\leq k}\\
a_{1}+a_{j-1}-a_{j-1}  & =-1,\;j\in\left\{  q_{i},2m+2-q_{i}\right\}  _{1\leq
i\leq k}\\
a_{j}+a_{2m+1-j}  & =a_{j+1}+a_{2m-j},\;2\leq j\leq m-1
\end{align}
and let $S^{\prime}\left(  m,q_{1},..q_{k}\right)  $ be the system given by
$\left(  1\right)  $ and $\left(  2\right)  $. We claim that any solution of
$S^{\prime}\left(  m,q_{1},..,q_{k}\right)  $ is also a solution of $S\left(
m,q_{1},..,q_{k}\right)  $. In fact, if $j\notin\left\{  q_{i},2m+2-q_{i}%
\right\}  _{1\leq i\leq k}$ then $2m+2-q_{i}\neq2m+1-j$ for all $i$, and
therefore we have
\[
a_{j}+a_{2m+1-j}=a_{j+1}+a_{2m-j}=a_{1}+a_{j}+a_{2m-j}%
\]
thus
\[
a_{2m+1-j}=a_{1}+a_{2m-j}%
\]
If $j\in\left\{  q_{i},2m+2-q_{i}\right\}  _{1\leq i\leq k}$ then
$2m+2-q_{i}=2m+1-j$ for some $i\in\left\{  1,..,k\right\}  $. Then
\[
a_{j}+a_{2m+1-j}=a_{j+1}+a_{2m-j}=1+a_{1}+a_{j}+a_{2m-j}%
\]
and simplifying
\[
a_{2m+1-j}=1+a_{1}+a_{2m-j}%
\]
This shows that the equations $\left(  3\right)  $ are superfluous. Now, the
system  $S^{\prime}\left(  m,q_{1},..q_{k}\right)  $ clearly has integer
solutions, depending on the parameters $a_{2}=N_{1}$ and $a_{3}=N_{2}$. Let
$\left(  a_{1},..,a_{2m}\right)  $ be the solution corresponding to the values
$N_{1}=N_{2}=1$ and $T_{2m+1}\left(  \mathbb{C}\left(  t\right)  \right)  $
denote the Borel subgroup of $GL\left(  2m+1,\mathbb{C}\right)  $ consisting
of lower triangular matrices. Define $f_{\left(  q_{1},..,q_{k}\right)  ,t}\in
T_{2m+1}\left(  \mathbb{C}\left(  t\right)  \right)  $ by
\[
\left\{
\begin{array}
[c]{l}%
f_{\left(  q_{1},..,q_{k}\right)  ,t}\left(  X_{i}\right)  =t^{a_{i}}%
X_{i},\;i\neq2m+1\\
f_{\left(  q_{1},..,q_{k}\right)  ,t}\left(  X_{2m+1}\right)  =t^{1+a_{2m-1}}%
\end{array}
\right.
\]
and consider the Lie algebra  $f_{\left(  q_{1},..,q_{k}\right)  ,t}^{-1}%
\ast\mu$, where $\mu$ is the Lie algebra law associated to $\frak{g}_{m}$.
Then the structural equations of  $f_{\left(  q_{1},..,q_{k}\right)  ,t}%
^{-1}\ast\mu$ are given by
\begin{align*}
d\omega_{1}  & =d\omega_{2}=0\\
d\omega_{j}  & =t^{a_{1}+a_{j-1}-a_{j}}\,\omega_{1}\wedge\omega_{j-1}%
,\;j\in\{q_{i},2m+2-q_{i}\}_{1\leq i\leq k}\\
d\omega_{j}  & =\omega_{1}\wedge\omega_{j-1},\;3\leq j\leq2m,\;j\notin
\{q_{i},2m+2-q_{i}\}_{1\leq i\leq k}\\
d\omega_{2m+1}  & =\sum_{j=2}^{m}\left(  -1\right)  ^{j}\omega_{j}\wedge
\omega_{2m+1-j}%
\end{align*}
Now, as $a_{1}+a_{j-1}-a_{j}=-1$ for $j\in\{q_{i},2m+2-q_{i}\}_{1\leq i\leq
k}$ by the system $S^{\prime}\left(  m,q_{1},..,q_{k}\right)  $ ( therefore by
$S\left(  m,q_{1},..,q_{k}\right)  $), it follows easily that
\[
\lim_{t\rightarrow\infty}\;f_{\left(  q_{1},..,q_{k}\right)  ,t}^{-1}\ast
\mu=\mu\left(  q_{1},..,q_{k}\right)
\]
where $\mu\left(  q_{1},..,q_{k}\right)  $ is the law associated to
$\frak{g}_{m}\left(  q_{1},..,q_{k}\right)  $.
\end{proof}

\begin{corollary}
For $m\geq4$%
\[
\mu\left(  q_{1},..,q_{k}\right)  \in O\left(  \mu\left(  q_{1}^{\prime
},..,q_{k}^{\prime}\right)  \right)
\]
if and only if $\left\{  q_{1},..,q_{k}\right\}  =\left\{  q_{1}^{\prime
},..,q_{k}^{\prime}\right\}  $.    
\end{corollary}

\begin{corollary}
For $m\geq 4$ the algebras $\frak{g}_{m}\left( q_{1},..,q_{k}\right) $ (
included $\frak{g}$) are nontrivial deformations of the algebra $\frak{h}%
_{m-1}\oplus \mathbb{C}^{2}$, where $\frak{h}_{m-1}$ is the $\left(
2m-1\right) $-dimensional Heisenberg Lie algebra.
\end{corollary}

\begin{proof}
Let $f_{t}\in T_{2m+1}\left( \mathbb{C}\left( t\right) \right) $ be defined
by 
\begin{equation*}
f_{t}\left( X_{i}\right) =t^{a_{i}}X_{i},\;1\leq i\leq 2m+1
\end{equation*}
where the $a_{i}$ satisfy the system $S^{\prime }\left(
m,q_{1},..,q_{m}\right) $ and the additional equation
\begin{equation*}
a_{1}+a_{2m-1}=a_{2m}-1
\end{equation*}
Then $\lim_{t\rightarrow \infty }f_{t}^{-1}\ast \mu \in O\left( \frak{h}%
_{m-1}\oplus \mathbb{C}^{2}\right) $.\newline A similar reasoning shows that $\frak{g%
}_{m}\left( q_{1},..,q_{k}\right) \rightarrowtail \frak{h}_{m-1}\oplus 
\mathbb{C}^{2}$. The result follows from the fact that a contraction defines
a nontrivial deformation [7].
\end{proof}

\section{Applications to complete Lie algebras}

Recall that a Lie algebra $\frak{g}$ is called complete if it is centerless
and any derivation is inner [8]. I recent years, a general theory of complete
solvable Lie algebras whose nilradical is of maximal rank has been developed
( [9], [10], [12]), and the existence of non-maximal rank algebras has been
pointed out. In this section we will see how to obtain families of
non-maximal rank completable Lie algebras considering the contractions above.%
\newline
Let $\frak{g}$ be a nilpotent Lie algebra and $Der\left( \frak{g}\right) $
its Lie algebra of derivations. A torus $\frak{t}$ over $\frak{g}$ is an
abelian subalgebra of $Der\left( \frak{g}\right) $ consisting of semisimple
derivations. The torus $\frak{t}$ induces a natural represntation on $\frak{g}
$ such that 
\begin{equation*}
\frak{g=}\sum_{\alpha \in \frak{t}^{\ast }}\frak{g}_{\alpha }
\end{equation*}
where $\frak{t}^{\ast }=Hom\left( \frak{t,}\mathbb{C}\right) $ and $\frak{g}%
_{a}=\left\{ X\in \frak{g\;}|\;[t,X]=\alpha \left( t\right) X,\;\forall t\in 
\frak{t}\right\} $. If $\frak{t}$ is maximal for the inclusion relation, by
the conjugation theorems of Morozov, its common dimension is a numerical
invariant called the rank of $\frak{g}$, denoted $rank\left( \frak{g}\right) 
$. An algebra is called of maximal rank if $rank\left( \frak{g}\right)
=b_{1}=\dim H^{1}\left( \frak{g},\mathbb{C}\right) $. Following Favre [6],
a weight system of $\frak{g}$ is given by 
\begin{equation*}
P\frak{g}\left( \frak{t}\right) =\left\{ \left( \alpha ,d\alpha \right)
\;|\;\alpha \in t^{\ast }\;\text{such that }\frak{g}_{\alpha }\neq
0,\;da=\dim \frak{g}_{\alpha }\right\} 
\end{equation*}
We also recall the folowing results:

\begin{proposition}
Let $\frak{g}$ be of maximal rank and $\frak{t}$ a maximal torus. Then
the semidirect product $\frak{t\oplus g}$ is complete solvable.
\end{proposition}

\begin{definition}
The algebra $\frak{g}$ is called completable if $\frak{t\oplus g}$ is
complete for a maximal torus $\frak{t}$, and simple completable if it is not the direct sum of nontrivial completable Lie algebras.
\end{definition}

The main result we will use is a slight modification of the next

\begin{theorem}
Let $\frak{g}$ be a Lie algebra and $\frak{h}$ a Cartan subalgebra. Assume
that following conditions are satisfied:

\begin{enumerate}
\item  $\frak{h}$ is abelian

\item  $\frak{g=}\sum_{\alpha \in \Delta }\frak{g}_{\alpha }$ with $\Delta
\subset \frak{h}^{\ast }-\{0\}$

\item  there is a basis $\left\{ \alpha _{1},..,\alpha _{l}\right\} $ of $%
\frak{h}^{\ast }$ n $\Delta $ such that $\dim \frak{g}_{\pm \alpha _{j}}\leq
1$ for $1\leq j\leq l$ and $[\frak{g}_{\alpha _{j}},\frak{g}_{-\alpha
_{j}}]\neq 0$ if $-\alpha _{j}\in \Delta $

\item  $\frak{h}$ and $\left\{ \frak{g}_{\pm \alpha _{j}},\;1\leq j\leq
l\right\} $ generate $\frak{g}$
\end{enumerate}
Then $\frak{g}$ is a complete Lie algebra.
\end{theorem}

The last conditions can be replaced by by more general statements [11] 

\begin{enumerate}
\item[3']  there is a enerating system $\left\{ \alpha _{1},..,\alpha
_{l}\right\} $ of $\frak{h}^{\ast }$ in $\Delta $ such that $\dim \frak{g}%
_{\alpha _{j}}=1$ for all $j$ and $\frak{h},\frak{g}_{\alpha _{1}},..,\frak{g%
}_{\alpha _{l}}$ genrate $\frak{g}$.

\item[4']  Let $0\neq x_{j}\in \frak{g}_{\alpha _{j}}$ and a basis $\left\{
\alpha _{1},..,\alpha _{r}\right\} $ of $\frak{h}^{\ast }$. For $r+1\leq
s\leq l,$%
\begin{equation*}
\alpha _{s}=\sum_{i=1}^{t}k_{is}\alpha _{j_{i}}-\sum_{i=1+t}^{r}k_{is}\alpha
_{j_{i}}
\end{equation*}
where $k_{is}\in \mathbb{N\cup \{}0\},\left( j_{1},..,j_{r}\right) $ is a
permutation  of $\left( 1,..,r\right) $, and there is a formula
\begin{eqnarray*}
&&\lbrack \underset{k_{1s}}{\underbrace{x_{j_{1}},..,x_{j_{1}}}},..\underset{%
k_{ts}}{\underbrace{x_{j_{t}},..,x_{j_{t}}},}...,x_{k_{m}}] \\
&=&[\underset{k_{t+1s}}{\underbrace{x_{j_{t}+1},..,x_{j_{t}+1}}},..\underset{%
k_{rs}}{\underbrace{x_{j_{r}},..,x_{j_{r}}},x_{s},x_{k_{1}}}...,x_{k_{m}}]
\end{eqnarray*}
without regard to the order or the way of bracketing, where $m\neq 0$ if $t=r
$.
\end{enumerate}

\begin{theorem1}
Let $\frak{g}$ be a Lie algebra satisfying conditions $\left( 1\right)
,\left( 2\right) ,\left( 3^{\prime }\right) $ and $\left( 4^{\prime }\right) 
$. Then $\frak{g}$ is a complete Lie algebra.
\end{theorem1}

We obtain one more consecuence of the theorem in the preceding section

\begin{corollary}
For $m\geq4$ and $k\geq1$
\[
rank\left(  \frak{g}_{m}\left(  q_{1},..,q_{k}\right)  \right)  >2
\]
\end{corollary}

\begin{proof}
The fact follows from the linear system $\left(  S\left(  \frak{g}_{m}\left(
q_{1},..,q_{k}\right)  \right)  \right)  $ associated to the algebras [1] and
the fact that $b_{1}\left(  \frak{g}_{m}\left(  q_{1},..,q_{k}\right)
\right)  >b_{1}\left(  \frak{g}_{m}\right)  =2$, where $b_{1}\left(
\frak{g}\right)  =\dim\;H^{1}\left(  \frak{g},\mathbb{C}\right)  $. 
\end{proof}

\begin{proposition}
Let $m\geq 4$ and $k\geq 1$ a weight system of $\frak{g}_{m}\left(
q_{1},..,q_{k}\right) $ is given by 
\begin{equation*}
P\frak{g}_{m}\left( \frak{t}_{m}\right) =\left\{ \left( \alpha _{i},d\alpha
_{i}\right) _{1\leq i\leq 2m+1}\right\} 
\end{equation*}
where $\dim \frak{g}_{\alpha _{i}}=1$ for all $i$ and the weights $\left\{
\alpha _{1},..,\alpha _{2m+1}\right\} $ satisfy the following linear system
\begin{align}
\alpha_{1}+\alpha_{j-1} &  =a_{j}, & \;3 &  \leq j\leq2m,\;j\notin\left\{
q_{i},2m+2-q_{i}\right\}  _{1\leq i\leq k}\nonumber\\
a_{2m-t}+-\alpha_{t+1}-\alpha_{m} &  =\alpha_{m+1}, & \;1 &  \leq t\leq
m-2\tag{$S_1$}%
\end{align}
In particular, $rank\left( \frak{g}_{m}\left( q_{1},..,q_{k}\right) \right)
\leq m+1$.
\end{proposition}

\begin{proof}
The system $\left( S_{1}\right) $ coincides with the linear system $S\left( 
\frak{g}_{m}\left( q_{1},..,q_{k}\right) \right) $ associated to the
nilpotent Lie algebra $\frak{g}_{m}\left( q_{1},..,q_{k}\right) $, thus the $%
\alpha _{i}$ correspond to eigenvalues of semisimple derivations [1].
Therefore we can isolate $\alpha _{j_{1}},..,\alpha _{j_{s}}$ $\left( 1\leq
j_{1}<j_{2}..<j_{s}<2+2k\right) $ such that for any $j\in \left\{
1,,,2m+1\right\} -\{j_{1},..,j_{s}\}$ we have 
\begin{equation*}
\alpha _{j}=\sum_{t=1}^{s}a_{j}^{t}\alpha _{j_{t}},\;a_{j}^{t}\in \mathbb{C}
\end{equation*}
Now, using the standard techniques [1] it is routine to verify that for $%
1\leq t\leq s$ the derivations $f_{j_{t}}\in Der\left( \frak{g}_{m}\left(
q_{1},..,q_{k}\right) \right) $ given by 
\begin{equation*}
f_{j_{t}}\left( X_{i}\right) =a_{i}^{t}X_{i},\;1\leq i\leq 2m+1
\end{equation*}
define a maximal torus $\frak{t}_{m}\left( q_{1},..,q_{k}\right) $ of $\frak{%
g}_{m}\left( q_{1},..,q_{k}\right) $. Thus the weight system is given as
above.\newline
Clearly any weight space $\frak{g}_{\alpha _{j}}$ is at most one
dimensional. For the last assertion, observe that for $k=m+1$ the system $%
\left( S_{1}\right) $ is 
\begin{eqnarray*}
\alpha _{1}+\alpha _{2m-1} &=&\alpha _{2m} \\
\alpha _{t}+\alpha _{2m+1-t} &=&\alpha _{m}+\alpha _{m+1},\;2\leq t\leq m-1
\end{eqnarray*}
and that the rank is the maximal possible, namely $m+1$.
\end{proof}

\begin{theorem}
For $m\geq 4$ and $k\geq 1$ the semidrect products 
\begin{equation*}
\frak{r}_{m}\left( q_{1},..,q_{k}\right) =\frak{t}_{m}\left(
q_{1},..,q_{k}\right) \oplus \frak{g}_{m}\left( q_{1},..,q_{k}\right) 
\end{equation*}
are solvable and complete.
\end{theorem}

\begin{proof}
It is easy to see that for any $f\in Der\left( \frak{r}_{m}\left(
q_{1},..,q_{k}\right) \right) $ we have 
\[
f\left( \frak{t}_{m}\left(
q_{1},..,q_{k}\right) \right) \subset \frak{g}_{m}\left(
q_{1},..,q_{k}\right) 
\]
This is a direct consequence of the particular
action of $\frak{t}_{m}\left( q_{1},..,q_{k}\right) $ over $\frak{g}%
_{m}\left( q_{1},..,q_{k}\right) $. If  $\frak{t}_{m}\left(
q_{1},..,q_{k}\right) =\left\{ h_{1},..,h_{s}\right\} $, then for any $1\leq
t\leq s$ there exists a permutation $\sigma \in \mathcal{S}_{s}$ such that $%
f_{j_{t}}=ad\left( h_{t}\right) $. The nilradical is clearly generated by $%
\left\{ X_{1},X_{2},X_{q_{i}},X_{2m+2-q_{i}}\right\} _{1\leq i\leq k}$,
where $\left\{ X_{1},..,X_{2m+1}\right\} $ is a dual basis to $\left\{
\omega _{1},..,\omega _{2m+1}\right\} $. It follows that the algebra $\frak{r%
}_{m}\left( q_{1},..,q_{k}\right) $ satisfies the conditions of the preceding theorem,
so that it is complete$.$
\end{proof}

\begin{lemma}
If $k=1$ and $q_{1}\neq m+1$ or $k\geq 2$, the $\frak{g}_{m}\left(
q_{1},..,q_{k}\right) $ is not of maximal rank.
\end{lemma}

\begin{proof}
The proof is an immediate consequence of the weight system. For $k=1$ and $%
q_{1}=m+1$ we have $rank(S_{1})=b_{1}=3$. Observe that this is the only 
case where $q_{i}=2m+2-q_{i}$.
\end{proof}

\begin{corollary}
For any odd dimension $n\geq 9$ there exist completable Lie algebras of
non-maximal rank.
\end{corollary}

Thus the contractions of the algebra $\frak{g}_{m}$ ( it is itself
completable since it is of maximal rank) are completable of non-maximal
rank, up to an exception. The families above expand the examples obtained in
[10],[12] for non-maximal rank. In fact, since the algebras $\frak{g}%
_{m}\left( q_{1},..,q_{k}\right) $ are nonsplit, we obtain even more:

\begin{corollary}
For $m\geq 4$ and $k\geq 1$ the algebras $\frak{g}_{m}\left(
q_{1},..,q_{k}\right) $ are simple completable.
\end{corollary}


\begin{thebibliography}{99}
\bibitem{}J. M. Ancochea, M. Goze. Le rang du syst\`{e}me lineaire des racines
d'une alg\`{e}bre de Lie resoluble complexe. Comm. Algebra \textbf{20 }(1992), 875-887.

\bibitem {}J. M. Ancochea, R. Campoamor. On Lie algebras whose nilradical is
(n-p)-filiform. Comm. Algebra \textbf{29} (2001), 427-450.

\bibitem {}J. M. Ancochea, R. Campoamor. On certain families of naturally
graded Lie algebras. J. Pure Appl. Algebra, \textit{to appear.}

\bibitem {}J. M. Ancochea, R. Campoamor. Nonfiliform characteristically
nilpotent and complete Lie algebras. Algebra. Colloq., \textit{to appear.}

\bibitem {}D. Burde. Degenerations of nilpotent Lie algebras. J. Lie Theory
\textbf{9} (1999), 193-202.

\bibitem {}G. Favre  Syst\`{e}me des poids sur une alg\`{e}bre de Lie nilpotente. Manuscripta Math.\textbf{9} (1973), 53-90.

\bibitem {}A. Fialowski, J. O'Halloran. A comparison of deformations and orbit closure. Comm. Algebra \textbf(18) (1990), 4121-4140.

\bibitem {}G. F. Leger.  Derivations of Lie algebras III. Duke Math. J.
\textbf{30} (1963), 637-645.

\bibitem {}D. J. Meng, L. S. Zhu. Solvable complete Lie algebras I. Comm.
Algebra \textbf{24} (1996), 4181-4197.

\bibitem {}D. J. Meng, L. S. Zhu. Solvable complete Lie algebras II. Algebra
Colloq. \textbf{5} (1998), 289-296.

\bibitem {}D. J. Meng. Complete Lie algebras and Heisenberg Lie algebras.
Comm. Algebra \textbf{22} (1994), 5509-5524.

\bibitem {}D. J. Meng. The complete Lie algebras with abelian nilpotent
radical. Acta Math. Sinica \textbf{34} (1991), 191-202.

\bibitem {}E. J. Saletan. Contractions of Lie groups. J. Math. Phys. \textbf{2
}(1961), 1-21. 

\bibitem {}C. Seeley. Degenerations of central quotients. Arch. Math. \textbf{56} (1991), 236-241.

\end{thebibliography}
\end{document}